\documentclass[10pt]{article}
\usepackage[utf8]{inputenc}

\usepackage{geometry,amsmath,mathtools,blkarray,amssymb,amsthm,setspace,xcolor,bbm,tikz-cd,xparse,amsfonts,authblk,graphics,datetime}

\usepackage{tikz,tikz-3dplot}
\usetikzlibrary{decorations.pathreplacing}

\usepackage[shortlabels]{enumitem}
\usepackage[colorlinks,linkcolor=blue,citecolor=blue,urlcolor=blue,bookmarks=false,hypertexnames=true]{hyperref}
\usepackage[backend=biber, style=bwl-FU, sortcites=false, maxcitenames=2, mincitenames=1, maxbibnames=8, uniquelist=false]{biblatex}
\usepackage{authblk}
\allowdisplaybreaks

\bibliography{smithson.bib}

\newtheorem{theorem}{Theorem}

\newtheorem{example}[theorem]{Example}
\newtheorem{remark}[theorem]{Remark}

\NewDocumentCommand{\defmathletter}{m}{%
    \expandafter\newcommand\csname b#1\endcsname{\mathbb{#1}}%
    \expandafter\newcommand\csname c#1\endcsname{\mathcal{#1}}%
}
\NewDocumentCommand{\defmathletters}{>{\SplitList{,}}m}{\ProcessList{#1}{\defmathletter}}
\defmathletters{A,B,C,D,E,F,G,H,I,J,K,L,M,N,O,P,Q,R,S,T,U,V,X,Y,Z}

\NewDocumentCommand{\defvector}{m}{%
    \expandafter\newcommand\csname v#1\endcsname{\mathbf{#1}}%
}
\NewDocumentCommand{\defvectors}{>{\SplitList{,}}m}{\ProcessList{#1}{\defvector}}
\defvectors{A,B,C,D,E,F,G,H,I,J,K,L,M,N,O,P,Q,R,S,T,U,V,X,Y,Z}
\defvectors{a,b,c,d,e,f,g,h,i,j,k,l,m,n,o,p,q,r,s,t,u,v,x,y,z}

\DeclareMathOperator{\gr}{graph}

\newdateformat{monthyeardate}{\monthname[\THEMONTH] \THEYEAR}

\title{On Smithson's fixed point theorem for order preserving multifunctions\footnote{Comments by Yifan Dai are gratefully acknowledged. This note supersedes an earlier draft by the first author.}}
\author[1]{Haruki Kono}
\affil[1]{\small Department of Economics, MIT, 77 Massachusetts Avenue, E52-402 Cambridge, MA 02139, USA, \url{hkono@mit.edu}}

\author[2]{Mark Voorneveld}
\affil[2]{\small{Department of Economics, Stockholm School of Economics, Box 6501, 11383 Stockholm, Sweden, \url{mark.voorneveld@hhs.se}}}
\date{\monthyeardate\today}

\begin{document}

\maketitle

\begin{abstract}
Fixed point theorems are ubiquitous in economic research.
Many studies cite Smithson (1971) ``Fixed points of order preserving multifunctions,'' yet the original proof contains errors.
This note presents a new, concise proof and explains why Smithson's argument is invalid. It also contains new results on the structure of the set of fixed points and monotone comparative statics.
\end{abstract}

\section{Introduction}

\cite{smithson1971} provides a fixed point theorem for multifunctions/correspondences on a partially ordered set. This theorem underlies many recent results on the existence and monotone comparative statics of, for instance, Nash equilibria in games with strategic complementarities \parencite{chekimkojima2021} and voting games \parencite{khan2020}, equilibria in large dynamic economies \parencite{acemoglujensen2015}, steady states of misspecified Bayesian learning models \parencite{ghosh2025}, and stable matchings in two-sided markets \parencite{chekimkojima2021}.

Unfortunately, its proof contains a mistake. After a brief reminder of notation and terminology, we therefore state the theorem, provide a new proof, and conclude by pointing out what goes wrong in Smithson's original argument. \cite{hoft1987} also suggests that the proof ``appears to be incorrect'' (p. 286), but does not establish this formally. She obtains Smithson's result as a corollary of her fixed point theorem for compositions of multifunctions. See Theorem 2.2 and Corollary 2.5 in \cite{hoft1987}. We provide a shorter, direct proof. Finally, we provide new results on the structure of the set of fixed points and monotone comparative statics.

\section{Notation}

A \textit{partially ordered set} is a set $X$ with a binary relation $\leq$ that is reflexive (for all $x \in X$, $x \leq x$), transitive (for all $x, y, z \in X$, $x \leq y$ and $y \leq z$ imply $x \leq z$), and antisymmetric (for all $x, y \in X$, $x \leq y$ and $y \leq x$ imply $x=y$). We write $x < y$ when $x \leq y$ and $x \neq y$. Let $Y$ be a subset of $X$. Element $y \in Y$ is \textit{maximal} in $Y$ if there is no $z \in Y$ with $y < z$. We call $x \in X$ an \textit{upper bound} on $Y$ if $y \leq x$ for all $y \in Y$ and a \textit{least upper bound} if, moreover, $x \leq z$ for all other upper bounds $z$. The least upper bound of $Y$ is also called its \textit{supremum} and denoted $\sup Y$. Set $Y$ is a \textit{chain} if its elements are totally ordered (for all $x, y \in Y$, $x \leq y$ or $y \leq x$).

The \textit{graph} of a function $f: X \to Y$ is the set $\gr(f) = \{(x, y) \in X \times Y \mid y = f(x)\}$. Function $f$ is \textit{isotone} if $X$ and $Y$ are partially ordered and for all $x, y \in X$, $x \leq y$ implies $f(x) \leq f(y)$.

A \textit{correspondence} or \textit{multifunction} $F : X \rightrightarrows X$ assigns a subset $F(x)$ of $X$ to each $x \in X$. A point $x \in X$ is a \textit{fixed point} of $F$ if $x \in F(x).$

\section{Smithson's fixed point theorem}

We now state the fixed point theorem in \cite{smithson1971} with our new proof.

\begin{theorem}\label{thm:smithson}
Let $(X, \leq)$ be a partially ordered set and $F: X \rightrightarrows X$ a multifunction. Assume
\begin{enumerate}[label=(A\arabic*),itemsep=.5pt,series=A]
\item Each nonempty chain in $X$ has a least upper bound; \label{ccpo}
\item There exist $e \in X$ and $y \in F(e)$ with $e\leq y$; \label{prefix}
\item For all $x_1,x_2\in X$ with $x_1\leq x_2$ and all $y_1\in F(x_1)$ there is a $y_2\in F(x_2)$ with $y_1\leq y_2$; \label{mono}
\item For each chain $C\subseteq X$ and each isotone function $f:C\to X$ with $f(x)\in F(x)$ for all $x\in C$, if $c =\sup C$, then there is a $y\in F(c)$ with $f(x)\leq y$ for all $x\in C$. \label{cont}
\end{enumerate}
Then $F$ has a fixed point, i.e., $x^* \in F(x^*)$ for some $x^* \in X$.
\end{theorem}
\begin{proof}
We start with a standard application of Zorn's lemma. Consider the set of functions
\[
\mathbf{F} \coloneqq \left\{ f: C \to X \ \ \middle| \ \
    \begin{aligned}
    & C \subseteq X, C \text{ is a nonempty chain}; \\
    & f \text{ is isotone}; \\
    & \text{for all } x \in C, \ f(x) \in F(x) \text{ and } x \leq f(x)
    \end{aligned}
\right\}.
\]
Set $\mathbf{F}$ is nonempty: given $e$ and $y$ as in \ref{prefix}, it contains the function $f: \{e\} \to X$ with $f(e) = y$. Partially order $\mathbf{F}$ by extension: for $f: C \to X$ and $g: C' \to X$ in $\mathbf{F}$, $f \leq g$ means $\gr(f) \subseteq \gr(g)$ (equivalently, $C \subseteq C'$ and $f = g$ on $C$). Each chain in $\mathbf{F}$ has an upper bound in $\mathbf{F}$ (e.g., the function whose graph is the union of those in the chain), so by Zorn's lemma $\mathbf{F}$ has a maximal element $f: C \to X$.

By \ref{ccpo}, $C$ has a least upper bound $x^*$. We will show that $x^*$ is a fixed point of $F$. By \ref{cont}, there is a $y \in F(x^*)$ with $f(x) \leq y$ for all $x \in C$. Since $x \leq f(x)$ for all $x \in C$ by the third condition of $\mathbf{F}$, $y$ is an upper bound on $C$. As $x^*$ is its least upper bound, $x^* \leq y$ holds.

If $x^* < y$, then we could extend $f$ to a larger element of $\mathbf{F}$, contradicting $f$'s maximality, as follows. Since $x \leq x^* < y$ for all $x \in C$, $C \cup \{y\}$ is a chain in $X$. By \ref{mono} (with $x_1 = x^*$ and $x_2 = y_1 = y$), there is a $y' \in F(y)$ with $y \leq y'$. Consequently, $x \leq f(x) \leq y \leq y'$ for all $x \in C$. So if we extend $f$ to $C \cup \{y\}$ by choosing $f(y) = y'$, it satisfies all conditions in $\mathbf{F}$: it is isotone and for all $x \in C \cup \{y\}$, $f(x) \in F(x)$ and $x \leq f(x)$. This extension contradicts that $f$ is maximal in $\mathbf{F}$. So it is not the case that $x^* < y$.

Since $x^* \leq y$ but not $x^* < y$, we get $x^* = y \in F(x^*)$, the desired fixed point.
\end{proof}

\begin{remark}
Condition \ref{cont} can be slightly weakened: we invoke it only for a function $f$ with the additional property that $x \le f(x)$ for all $x$ in its domain.
\end{remark}

Next, we discuss why the original proof of \cite{smithson1971} is incorrect. He argues that the collection
\begin{align*}
    \mathcal{S}
    \coloneqq
    \left\{
        Y \subseteq X
        \mid
        (1) \ e \in Y; \
        (2) \ x \in Y, \ e \leq z \leq x \Rightarrow z \in Y; \
        (3) \ x \in Y \Rightarrow \exists z \in F(x) \text { s.t. } x \leq z
    \right\}
\end{align*}
of subsets of $X$, partially ordered by set inclusion, has a maximal element $X_0$ and that there is a maximal chain $C$ in $X_0$. In the third paragraph of his proof he claims that, regardless of the choice of $X_0$ and $C$, the supremum of $C$ lies in $X_0$. This claim, crucial to his remaining steps, is wrong.

To give a counterexample, let $X \coloneqq [0,1] \cup \{a\}$ be the real interval $[0,1]$ with an extraneous element $a$. Partially order $X$ by ordering elements of $[0,1]$ as usual and declaring that $0 < a < 1$, but $a$ is incomparable to any other element in $(0,1)$. Formally,
\[
x \leq y \qquad \iff \qquad
\begin{cases}
x, y \in [0, 1] \text{ and } x \leq y \text{ (in the usual order) or} \\
(x,y) \in \{ (0,a), \ (a,1), \ (a,a)\}
\end{cases}
\]
Define correspondence $F$ by $F(x) = \{x\}$ if $x \in [0,1]$ and $F(a)=\{0\}$.

All assumptions in the fixed point theorem are satisfied. Throughout this example we choose $e = 0$ to be the element referred to in \ref{prefix}. For \ref{cont}, note that a chain in $X$ cannot contain elements of $(0,1)$ as well as $a$; it is a subset of $[0,1]$ or of $\{0,a,1\}$.

Smithson's claim (if $X_0$ is a maximal element of $\mathcal{S}$ and $C$ is a maximal chain in $X_0$, then $C$'s supremum lies in $X_0$) is violated in this example. To see this, set $X_0 = [0,1).$
We observe that $X_0$ is a maximal element of $\mathcal{S}$. It is easy to confirm that it lies in $\mathcal{S}$. For maximality, note that the only elements of $X$ not in $X_0$ are $1$ and $a$. We cannot add $a$, as $a$ violates the third condition in $\mathcal{S}$. Likewise, we cannot add 1: the second condition in $\mathcal{S}$ then implies that we must add $a$ as well.
Since $X_0 = [0,1)$ is a chain, it is the unique maximal chain $C$ in $X_0$. But its least upper bound 1 does not lie in $X_0$, contradicting Smithson's claim.

\section{The structure of the set of fixed points and monotone comparative statics}

In this section, we generalize the results of \cite{chekimkojima2021} about the structure of the set of fixed points and monotone comparative statics in a number of directions. More precisely, our Theorems \ref{thm: structure} and \ref{thm: mcs} generalize Theorems 6 and 7 of \cite{chekimkojima2021}, respectively. First, our theorems hold under the conditions of Smithson's fixed point theorem, whereas \cite{chekimkojima2021} impose additional topological assumptions on $X$. Moreover, we not only establish the existence of maximal fixed points but show that there is one above every element $e$ satisfying \ref{prefix}. Dispensing with the need for topological tools, our proofs are considerably less involved.

First, there is a maximal fixed point above each element satisfying condition \ref{prefix}:

\begin{theorem}\label{thm: structure}
Under the conditions of Smithson's fixed point theorem: for each $e$ as in \ref{prefix}, the multifunction has a maximal fixed point $x^*$ with $e \le x^*$.
\end{theorem}
\begin{proof}
Our proof of Theorem \ref{thm:smithson} carries over verbatim if we add to $\mathbf{F}$ the condition that an element $e$ as in \ref{prefix} must lie in $C$. Given a maximal $f:C \to X$ in $\mathbf{F}$, that proof establishes that $x^* = \sup C$ is a fixed point and that $f(x) \le x^*$ for all $x \in C$. Since $e \in C$, $e \le \sup C = x^*$.

And $x^*$ is a maximal fixed point. If, to the contrary, $x^* < y^*$ for some fixed point $y^*$, then we could extend $f$ to $C \cup \{y^*\}$ by defining $f(y^*) = y^*$, contradicting $f$'s maximality in $\mathbf{F}$. The other properties of $\mathbf{F}$ are easy to verify, so we only argue why this extended $f$ remains isotonic. The added $y^*$ is the greatest element of the chain $C \cup \{y^*\}$: for each $x \in C$, $x \le \sup C = x^* < y^*$. And $f(y^*) = y^*$ is the greatest function value: recalling that $f(x) \le x^*$ for all $x \in C$, it follows that $f(x) \le x^* < y^* = f(y^*)$. This extension contradicts $f$'s maximality. Hence $x^*$ is a maximal fixed point.
\end{proof}

Theorem \ref{thm: structure} ensures the existence of maximal fixed points. The following examples show that minimal fixed points need not exist. Nor least or greatest ones; in fact, the set of fixed points may consist entirely of incomparable elements.

\begin{example}\label{ex: no min}
Multifunction $F: (0,1] \rightrightarrows (0,1]$ on the real interval $(0,1]$ with its usual order defined by $F(x) = \{x\}$ satisfies all conditions in Smithson's fixed point theorem. Its set of fixed points $(0,1]$ has no minimal element.
\end{example}

\begin{example}
Let $X$ consist of three elements, denoted $0$, $0'$, and $1$, partially ordered such that $0 < 1$, $0' < 1$, but $0$ and $0'$ are incomparable. Correspondence $F: X \rightrightarrows X$ with $F(0) = \{0\}$, $F(0') = \{0'\}$, and $F(1) = \{0,0'\}$ satisfies all conditions in Smithson's fixed point theorem. Its set of fixed points consists of the incomparable elements $0$ and $0'$.
\end{example}

Next, we formulate a monotone comparative statics result for fixed points. Informally, it claims that if multifunction $F$ is obtained from multifunction $G$ by ``moving it upwards'' (for each $x$ in the domain, if we have a point in $G(x)$, there is a larger point in $F(x)$), then the set of fixed points moves in the same direction.

\begin{theorem}\label{thm: mcs}
Let $(X, \le)$ be a partially ordered set and $F: X \rightrightarrows X$ and $G: X \rightrightarrows X$ two multifunctions. Assume that $F$ satisfies all conditions in Smithson's fixed point theorem and
\begin{enumerate}[resume*=A]
\item For each $x \in X$ and each $y \in G(x)$ there is a $z \in F(x)$ with $y \le z$. \label{mcs}
\end{enumerate}
Then for each fixed point $x_0$ of $G$ there is a maximal fixed point $z_0$ of $F$ with $x_0 \le z_0$.
\end{theorem}
\begin{proof}
Let $x_0$ be a fixed point of $G$: $x_0 \in G(x_0)$. By \ref{mcs} there is a $y_0 \in F(x_0)$ with $x_0 \le y_0$: $x_0$ and $y_0$ satisfy condition \ref{prefix} in Smithson's fixed point theorem. By Theorem \ref{thm: structure}, $F$ has a maximal fixed point $z_0$ with $x_0 \le z_0$.
\end{proof}

As a final remark, \cite{chekimkojima2021} refer to monotonicity condition \ref{mono} as ``weak upper monotonicity.'' It helps to ensure the existence of maximal fixed points, not necessarily minimal ones (see Example \ref{ex: no min}). Of course, dual results about minimal (not necessarily maximal) fixed points can be obtained by adjusting the argument to ``weak lower monotonicity'': For all $x_1, x_2 \in X$ with $x_1 \le x_2$ and all $y_2 \in F(x_2)$ there is a $y_1 \in F(x_1)$ with $y_1 \le y_2$.

\printbibliography

\end{document}